\documentclass[a4paper]{amsart}
\usepackage{graphicx}
\usepackage{amssymb}
\usepackage{amsmath}
\usepackage{amsthm}
\usepackage{amscd}
\usepackage[all,2cell]{xy}

\UseAllTwocells \SilentMatrices
\newtheorem{thm}{Theorem}[section]

\newtheorem{cor}[thm]{Corollary}
\newtheorem{lem}[thm]{Lemma}

\newtheorem{prop}[thm]{Proposition}
\theoremstyle{definition}

\theoremstyle{remark}
\newtheorem{rem}[thm]{\bf Remark}
\numberwithin{equation}{section}

\begin{document}
\title[The Stable Monomorphism category of a Frobenius category]
{The stable monomorphism category of a Frobenius category}
\author[  Xiao-Wu Chen
] {Xiao-Wu Chen}
\date{Dec. 10, 2009}
\thanks{This project was supported by Alexander von Humboldt
Stiftung. The author is also supported by China Postdoctoral Science
Foundation (No.s 20070420125 and 200801230) and by National Natural
Science Foundation of China (No.10971206).}

\thanks{E-mail:
xwchen$\symbol{64}$mail.ustc.edu.cn}

\maketitle

\dedicatory{}%
\commby{}%

\begin{center}
\end{center}

\begin{abstract}
For a Frobenius abelian category  $\mathcal{A}$, we show that
the category ${\rm Mon}(\mathcal{A})$ of  monomorphisms in $\mathcal{A}$ is a Frobenius exact category;
 the associated stable category $\underline{\rm Mon}(\mathcal{A})$ modulo projective objects is called the
 stable monomorphism category of $\mathcal{A}$. We show that a tilting object in the stable category
 $\underline{\mathcal{A}}$ of $\mathcal{A}$ modulo projective objects induces naturally
 a tilting object in $\underline{{\rm Mon}}(\mathcal{A})$. We show that if  $\mathcal{A}$ is the category
 of (graded) modules over a (graded) self-injective algebra $A$, then the stable monomorphism category
 is triangle equivalent to the (graded) singularity category of the (graded) $2\times 2$
 upper triangular matrix algebra $T_2(A)$. As an application, we give
two characterizations to the stable category of Ringel-Schmidmeier
(\cite{RS3}).
\end{abstract}

\section{Introduction}

Let $\mathcal{A}$ be an abelian category. Denote by  ${\rm Mor}(\mathcal{A})$
 the category of morphisms in $\mathcal{A}$ (\cite[p.101]{ARS}):
 the objects are morphisms in $\mathcal{A}$ and the morphisms are given  by commutative squares in $\mathcal{A}$.
 It is an abelian category (\cite[Proposition 1.1]{FGR}).
 We are mainly concerned with the full subcategory ${\rm Mon}(\mathcal{A})$ of  ${\rm Mor}(\mathcal{A})$
 consisting of monomorphisms in $\mathcal{A}$,
which is called the \emph{monomorphism category} of $\mathcal{A}$. It is  an additive
subcategory of ${\rm Mor}(\mathcal{A})$ which is  closed under extensions, thus it becomes
an exact category in the sense of Quillen (\cite[Appendix A]{Ke3}).

\vskip 5pt

   In the case that the abelian category $\mathcal{A}$ is the
module category over a ring, the monomorphism category ${\rm
Mon}(\mathcal{A})$ is known as the \emph{submodule category}.
Recently it is studied intensively by Ringel and  Schmidmeier
(\cite{RS1, RS2, RS3}). If the ring is the truncated polynomial ring
$\mathbb{Z}[t]/(t^p)$ with $p\geq 2$ and $t$ an indeterminant, the
study of the submodule category goes back to Birkhoff (\cite{B}; see
also \cite{Ar}). The case that the ring is $k[t]/(t^p)$ with $k$ a
field is studied by Simson (\cite{Sim}) and also by Ringel and
Schmidmeier. In this case, the study of indecomposable objects in
${\rm Mon}(\mathcal{A})$ shows
 an example of the typical trichotomy  phenomenon  ``finite/tame/wild" in the representation
 theory of finite dimensional algebras, where the trichotomy depends on the
parameter  $p$; see \cite[section 6]{RS3}. Moreover, the case where
the abelian category $\mathcal{A}$ is given by the graded module
category over  the graded algebra $k[t]/(t^p)$ with ${\rm deg}\;
t=1$ plays an important role in \cite{RS3}; in this case, the
monomorphism category ${\rm Mon}(\mathcal{A})$ is denoted by
$\mathcal{S}(\widetilde{p})$. It is a Frobenius exact category
(\cite{KML}; also see Lemma \ref{lem:Frobeniusexact} and compare
\cite[section 5]{Ke3}). Then by \cite[Chapter I, Theorem 2.6]{H1}
its stable category $\underline{\mathcal{S}}(\widetilde{p})$ modulo
projective objects is triangulated. A very recent and remarkable
result due to Kussin, Lenzing and Meltzer claims
 that the stable category  $\underline{\mathcal{S}}(\widetilde{p})$ is triangle equivalent
to the stable category of vector bundles on the weighted projective lines of type $(2, 3, p)$; see \cite{KML}. Note that
a similar trichotomy phenomenon  ``domestic/tubular/wild" occurs in
the classification of indecomposable vector bundles on the weighted projective lines of type $(2,3, p)$,
while the trichotomy again depends on the parameter $p$; see \cite{Len, KMLP}.
In this paper,  we will call the triangulated category $\underline{\mathcal{S}}(\widetilde{p})$
the \emph{stable category of Ringel-Schmidmeier}.

\vskip 5pt

 The present paper studies the monomorphism category ${\rm Mon}(\mathcal{A})$ of a Frobenius
 abelian category $\mathcal{A}$, in particular, the stable category $\underline{\mathcal{A}}$
 modulo projective objects
 is triangulated. We show that ${\rm Mon}(\mathcal{A})$ is a Frobenius exact category and then
 the stable category $\underline{\rm Mon}(\mathcal{A})$  modulo projective objects  is triangulated; it is called
 the \emph{stable monomorphism category} of $\mathcal{A}$. Recently this category
 is also studied by \cite{IKM}. Note that the triangulated categories above
 are algebraical in the sense of Keller.
We have a well-behaved notion
 of \emph{tilting object} for an algebraical triangulated category (\cite{Ke4}). We prove that a tilting object in  $\underline{\mathcal{A}}$ induces naturally a tilting object in $\underline{\rm Mon}(\mathcal{A})$; see Theorem \ref{thm:first}. Moreover, if the category $\mathcal{A}$ is
 the (graded) module category over a (graded) self-injective algebra $A$,
 we relate the category $\underline{\rm Mon}(\mathcal{A})$ to  the category of (graded)
 Gorenstein-projective modules and then to the (graded) singularity category of the $2\times 2$
 upper triangular matrix algebra $T_2(A)$ of $A$ (for $T_2(A)$, see \cite[p.115]{FGR} and
  \cite[Chapter III, section 2]{ARS}); see Theorem \ref{thm:second}. We are inspired by a
  computational result by Li and Zhang  on Gorenstein-projective modules (\cite{LZ}).
  Here the Gorenstein-projective module is in the sense of Enochs and Jenda (\cite[Chapter 10]{EJ}),
  and the singularity category is in the
  sense of Orlov (\cite{O1, O2}; compare \cite{Buc, H2}).

   \vskip 5pt

   Combining all these together,  we give two characterizations to the stable category of
    Ringel-Schmidmeier in Theorem \ref{thm:last}: we  characterize the stable category
  $\underline{\mathcal{S}}(\widetilde{p})$ as the bounded derived
  category of $T_2(k\mathbb{A}_{p-1})$, where $\mathbb{A}_{p-1}$ is the linear quiver
   with $p-1$ vertices and linear orientation, and $k\mathbb{A}_{p-1}$ is the path algebra;
   we  characterize  the stable category $\underline{\mathcal{S}}(\widetilde{p})$ as the graded singularity
   category of $T_2(k)[t]/(t^p)$, where the algebra $T_2(k)[t]/(t^p)$ is graded
   such that ${\rm deg}\; T_2(k)=0$ and ${\rm deg}\; t=1$.

\vskip 5pt

For the convention, throughout we fix a commutative artinian ring
$R$. All artin algebras
 are artin $R$-algebras, and all categories and functors are $R$-linear.
 For an artin algebra $A$, denote by
 ${\rm mod}\; A$ the category of finitely generated right $A$-modules and by ${\rm proj}\; A$  the
 full subcategory consisting of projective modules. We denote by $A_A$ and $_AA$ the right and
 left regular modules of the artin algebra $A$, respectively. For triangulated categories and derived categories,
 we refer to \cite{Har, H1, Ke3.5, Ke4}.

\section{Monomorphism Category}

Let $\mathcal{A}$ be a Frobenius abelian category. Thus
$\mathcal{A}$ has enough projective objects and enough injective
objects, and the class of projective objects coincides with the
class of injective objects. Denote by $\mathcal{P}$ the full
subcategory of $\mathcal{A}$ consisting of projective objects.
Denote by $\underline{\mathcal{A}}$ the \emph{stable category} of
$\mathcal{A}$ modulo $\mathcal{P}$: the objects are the same as
$\mathcal{A}$, and the morphism spaces are factors of the morphism
spaces in $\mathcal{A}$ modulo
 those factoring through projective objects (\cite[p.101]{ARS}).
The stable category $\underline{\mathcal{A}}$ is a triangulated
category such that its shift functor is given by the inverse of the
syzygy functor on $\underline{\mathcal{A}}$ and triangles are
induced by short exact sequences in $\mathcal{A}$; for details, see
\cite[Chapter I, section 2]{H1}.

\vskip 5pt

Recall that ${\rm Mor}(\mathcal{A})$ is the category of morphisms in $\mathcal{A}$: the objects are morphisms $\alpha \colon A\rightarrow B$ in $\mathcal{A}$ and the morphisms are commutative squares in $\mathcal{A}$, that is, of the form $(f, g)\colon \alpha\rightarrow \alpha'$ where $f\colon A \rightarrow A'$ and $g\colon B\rightarrow B'$ are morphisms in $\mathcal{A}$ such that $\alpha'\circ f=g\circ \alpha$ (\cite[p.101]{ARS}). For an object $\alpha\colon A\rightarrow B$ in  ${\rm Mor}(\mathcal{A})$,
 we write $s(\alpha)=A$ and $t(\alpha)=B$, which are called the \emph{source} and \emph{target} of $\alpha$, respectively. Note that ${\rm Mor}(\mathcal{A})$ is an abelian category such that a sequence $\alpha'\rightarrow \alpha\rightarrow \alpha''$ is exact if and only if the induced sequences of sources and targets are exact in $\mathcal{A}$ (\cite[Corollary 1.2]{FGR}).

\vskip 5pt

Recall that an exact category in the sense of Quillen is an additive
category together with an \emph{exact structure}, that is, a
distinguished class of  ker-coker sequences, which are called
\emph{conflations}, subject to certain axioms. Note that a full
additive subcategory of an abelian category which is closed under
extensions has a natural exact structure such that conflations are
just short exact sequences with terms in the subcategory
(\cite[Appendix A]{Ke3} and \cite[section 4]{Ke3.5}). Moreover,
there is a notion of Frobenius exact category and the associated
stable category modulo projective objects is still triangulated;
compare \cite[p.10-11]{H1},
 \cite[subsection 1.2 b)]{Ke3} and \cite[section 6]{Ke3.5}.

\vskip 5pt

 Recall that our main concern is  the \emph{monomorphism category} ${\rm Mon}(\mathcal{A})$,
 which is the full subcategory of ${\rm Mor}(\mathcal{A})$ consisting of monomorphisms in $\mathcal{A}$.
 We will consider the following two functors: the first
functor $i_1\colon \mathcal{A}\rightarrow {\rm Mon}(\mathcal{A})$ is
defined such that $i_1(A)=0\rightarrow A$ and $i_1(f)=(0, f)$ where
$A$ is an object and $f$ is a morphism in $\mathcal{A}$; the second
$i_2\colon \mathcal{A}\rightarrow {\rm Mon}(\mathcal{A})$ is defined
such that $i_2(A)=A\stackrel{{\rm Id}_A}\rightarrow A$ and
$i_2(f)=(f, f)$. Note that both functors are exact and fully
faithful.

\begin{lem}\label{lem:Frobeniusexact}
Let $\mathcal{A}$ be an abelian category.
Then the monomorphism  category ${\rm Mon}(\mathcal{A})$ is an exact category such that
conflations are given by sequences $\alpha' \rightarrow \alpha
\rightarrow \alpha''$ with the induced sequences of sources and
targets short exact in $\mathcal{A}$. \par
  Assume further that $\mathcal{A}$ is Frobenius. Then the exact category  ${\rm Mon}(\mathcal{A})$
is Frobenius such that its projective objects are equal to direct summands of objects of
the form $i_1(P)\oplus i_2(P)$ where $P$ is a projective object in
$\mathcal{A}$.
\end{lem}

\begin{proof}
Note that ${\rm Mon}(\mathcal{A})$ is an  additive subcategory of
the abelian category ${\rm Mor}(\mathcal{A})$ which is closed under
extensions by Snake Lemma. Then it is an exact category with
conflations induced by short exact sequences in  ${\rm
Mor}(\mathcal{A})$; see Example 4.1 in \cite{Ke3.5}.

\vskip 3pt

Assume now that the abelian category $\mathcal{A}$ is Frobenius.
We will show first that objects of the form $i_1(P)$ and $i_2(P)$
 are projective and injective.  Recall that for
 an object $\alpha$ in ${\rm Mon} (\mathcal{A})$ we denote
 by $s(\alpha)$ and $t(\alpha)$ the source and target of $\alpha$,
  respectively. Note that
$${\rm Hom}_{{\rm Mon}(\mathcal{A})}(i_1(P), \alpha)\simeq {\rm
Hom}_\mathcal{A}(P, t(\alpha)) \mbox{  and }{\rm Hom}_{{\rm
Mon}(\mathcal{A})}(i_2(P), \alpha)\simeq {\rm Hom}_\mathcal{A}(P,
s(\alpha)).$$ These isomorphisms show that the objects are
projective. Note that $${\rm Hom}_{{\rm Mon}(\mathcal{A})}( \alpha,
i_1(P))\simeq {\rm Hom}_\mathcal{A}({\rm Cok}\; \alpha, P) \mbox{
and } {\rm Hom}_{{\rm Mon}(\mathcal{A})}( \alpha, i_2(P))\simeq {\rm
Hom}_\mathcal{A}(t(\alpha), P).$$ These isomorphisms show that the
objects are injective (note that the functor ${\rm Cok}$ of taking
the cokernels is exact on ${\rm Mon}(\mathcal{A})$ by Snake Lemma).

\vskip 3pt

 Let $\alpha$ be an object in ${\rm Mon}(\mathcal{A})$. Take
epimorphisms $P\rightarrow s(\alpha)$ and $P\rightarrow t(\alpha)$
with $P$ projective in $\mathcal{A}$. Then we have an epimorphism
$i_1(P)\oplus i_2(P)\rightarrow \alpha$ whose kernel lies in ${\rm
Mon}(\mathcal{A})$. This shows that the exact category ${\rm Mon}
(\mathcal{A})$ has enough projective objects. On the other hand, for
the object $\alpha$, take monomorphisms  $a\colon
t(\alpha)\rightarrow P$ and $b'\colon {\rm Cok} \; \alpha
\rightarrow P$ with $P$ projective in $\mathcal{A}$. Denote by $b$
the composite $t(\alpha)\rightarrow {\rm Cok}\; \alpha
\stackrel{b'}\rightarrow P$ where the first morphism is the natural
projection. Consider the following morphism in ${\rm
Mor}(\mathcal{A})$
$$(\begin{pmatrix}a\circ \alpha \\ 0 \end{pmatrix}, \begin{pmatrix} a \\ b \end{pmatrix})\colon
\alpha\longrightarrow i_2(P)\oplus i_1(P).$$ It is a monomorphism
and by a diagram-chasing its cokernel lies in ${\rm
Mon}(\mathcal{A})$. Then it becomes a conflation in ${\rm
Mon}(\mathcal{A})$. This shows that the exact category ${\rm
Mon}(\mathcal{A})$ has enough injective objects. From the argument
above, it is direct to conclude that in the exact category ${\rm
Mon} (\mathcal{A})$ the class of projective objects coincides with
the class of injective objects, and projective objects are direct
summands of objects
 of the form $i_1(P)\oplus i_2(P)$ where $P$ is a projective
object in $\mathcal{A}$.
\end{proof}

\begin{rem}
With a slightly modified  proof as above, one can show that a
similar result holds if the category $\mathcal{A}$ is an exact
category. In this case, one replaces ${\rm Mon}(\mathcal{A})$ by the
category of \emph{inflations} in $\mathcal{A}$; compare
\cite[section 5]{Ke3} and \cite{IKM}. \hfill $\square$
\end{rem}

\vskip 5pt

For a Frobenius abelian category $\mathcal{A}$, we denote by ${\underline {\rm Mon}}(\mathcal{A})$
the stable category of ${\rm Mon}(\mathcal{A})$ modulo projective objects; it is a
triangulated category. We will call it the \emph{stable monomorphism
category} of $\mathcal{A}$.

\vskip 5pt

 Note that both the functors $i_1$ and $i_2$  are fully faithful and send projective objects
to projective objects. Then they induce fully faithful triangle
functors $i_1\colon \underline{\mathcal{A}}\rightarrow {\underline {\rm Mon}}(\mathcal{A})$ and
$i_2\colon \underline{\mathcal{A}}\rightarrow {\underline {\rm Mon}}(\mathcal{A})$ (\cite[p.23, Lemma 2.8]{H1}).

\section{Tilting objects in Stable Monomorphism Category}

In this section, we will show that for a Frobenius abelian category
$\mathcal{A}$, a tilting object in the stable category
$\underline{\mathcal{A}}$ induces naturally a tilting object in the
stable monomorphism category  ${\underline {\rm Mon}}(\mathcal{A})$.

\vskip 5pt

Following Keller we call that a triangulated category is
\emph{algebraical} provided that it is triangle equivalent to the
stable category of a Frobenius exact category (\cite[subsection
8.7]{Ke4}). One has a well-behaved notion of tilting object in an
algebraical triangulated category. \vskip 5pt

Let $\mathcal{T}$ be an algebraical triangulated category. Denote by
$[1]$ the shift functor and by $[n]$ its $n$-th power for each $n\in
\mathbb{Z}$. An object $T$ in $\mathcal{T}$ is a \emph{tilting
object}  if the following conditions are satisfied:

\begin{enumerate}
\item[(T1)] ${\rm Hom}_\mathcal{T}(T, T[n])=0$ for $n\neq 0$;
\item[(T2)] the smallest \emph{thick} triangulated subcategory of $\mathcal{T}$
containing $T$ is $\mathcal{T}$ itself;
\item[(T3)] ${\rm End}_\mathcal{T}(T)$ is an artin algebra
having finite global dimension.
\end{enumerate}
Here we recall that a triangulated subcategory of $\mathcal{T}$
is called \emph{thick} if it is closed under taking direct summands.
Note that the notion of tilting object presented here is slightly different from, however closely
related to, the ones in \cite{H1} and \cite{Ke4}.

\vskip 5pt

Recall that an additive category is said to be \emph{idempotent-split} provided that
 each idempotent $e\colon X\rightarrow X$ admits a factorization
 $X\stackrel{u}\rightarrow Y\stackrel{v}\rightarrow X$ such that $u\circ v={\rm Id}_Y$ (\cite[Chapter I, 3.2]{H1}).
Note that for an artin algebra $A$ having finite global dimension,
the bounded derived category $\mathbf{D}^b({\rm mod}\; A)$ is
algebraical and idempotent-split (see the proof of \cite[Chapter I,
Corollary 4.9]{H1}), and it has $A_A$ as its tilting object.

\vskip 5pt

 The following remarkable
result due to Keller claims that the converse holds true (compare
\cite[Theorem 1]{BK}).

\begin{lem} \label{lem:Keller} {\rm (Keller)} Let $\mathcal{T}$
be an idempotent-split algebraical triangulated
category with a tilting object $T$. Then there is a triangle
equivalence
$$\mathcal{T}\; \simeq \; \mathbf{D}^b({\rm mod}\; {\rm
End}_\mathcal{T}(T)).$$
\end{lem}

\begin{proof}
Set $A={\rm End}_\mathcal{T}(T)$. By \cite[Theorem 8.51 a)]{Ke4}
there is a triangle functor $F'\colon \mathcal{T}\rightarrow {\bf
D}(A')$ sending $T$ to $A'$, where $A'$ is a differential graded
algebra with the only nonzero cohomology $H^0(A')\simeq A$ and
$\mathbf{D}(A')$ is the (unbounded) derived category of differential
graded (right) modules on $A'$.  Note that by \cite[subsection
8.4]{Ke4} there is a triangle equivalence $\mathbf{D}(A')\simeq
\mathbf{D}(\mbox{Mod}\; A)$ identifying $A'$ with $A_A$, where ${\rm
Mod}\; A$ is the category of (unnecessarily finitely generated)
right $A$-modules. Consequently, there is a  triangle functor
$F\colon \mathcal{T}\rightarrow \mathbf{D}(\mbox{Mod}\; A)$ sending
$T$ to $A$. Using (T1) and (T2) and applying Beilinson Lemma
(\cite[p.72, Lemma 3.4]{H1}), the triangle functor $F$ is fully
faithful. Then we may view $\mathcal{T}$ as a triangulated
subcategory of $\mathbf{D}(\mbox{Mod}\; A)$; moreover, since
$\mathcal{T}$ is idempotent-split, it is necessarily a thick
subcategory of $\mathbf{D}(\mbox{Mod}\; A)$. By (T3) the artin
algebra $A$ has finite global dimension, and then the smallest thick
triangulated subcategory of  $\mathbf{D}(\mbox{Mod}\; A)$ containing
$A_A$ is $\mathbf{D}^b(\mbox{mod}\; A)$. From this we conclude that
the essential image of $F$ is $\mathbf{D}^b(\mbox{mod}\; A)$.
Therefore $F$ induces the required equivalence.
\end{proof}

Our first main observation states that a tilting object in the
stable category $\underline{\mathcal{A}}$ induces naturally a
tilting object in the stable monomorphism category $\underline{\rm
Mon}(\mathcal{A})$. Recall that for an artin algebra $A$,
$T_2(A)=\begin{pmatrix} A & A \\ 0 & A
\end{pmatrix}$ is the $2\times 2$ upper triangular  matrix algebra
(\cite[Chapter III, section 2]{ARS}).

\begin{thm}\label{thm:first}
Let $\mathcal{A}$ be a Frobenius abelian category such that $T$ is a
tilting object in its stable category $\underline{\mathcal{A}}$.
Then $T'=i_1(T)\oplus i_2(T)$ is a tilting object in $\underline{\rm
Mon}(\mathcal{A})$; moreover, we have an isomorphism ${\rm End}_{\underline{\rm
Mon}(\mathcal{A})}(T')\simeq T_2({\rm End}_{\underline{\mathcal{A}}}(T))$ of algebras.
\end{thm}

\begin{proof}
Recall that $i_1\colon \underline{\mathcal{A}}\rightarrow
\underline{\rm Mon}(\mathcal{A})$ and $i_2\colon
\underline{\mathcal{A}}\rightarrow \underline{\rm Mon}(\mathcal{A})$
are fully faithful triangle functors. Note that for objects $A$ and
$B$ in $\mathcal{A}$, ${\rm Hom}_{\underline{\rm Mon}(\mathcal{A})}
(i_2(A), i_1(B))=0$. So to check the condition (T1) for $T'$, it
suffices to show that ${\rm Hom}_{\underline{\rm Mon}(\mathcal{A})}
(i_1(T), i_2(T)[n])=0$ for $n\neq 0$. For this end, note that since
$i_2$ is a triangle functor, we have
$$i_2(T)[n] \simeq i_2 (T[n])=\; T[n]\stackrel{{\rm
Id}_{T[n]}}\longrightarrow T[n].$$ Thus a morphism in ${\rm
Hom}_{\underline{\rm Mon}(\mathcal{A})} (i_1(T), i_2(T)[n])$  is of
the form $(0, f)$, where $f\colon T\rightarrow T[n]$ is a morphism
in $\mathcal{A}$. By the condition (T1) for $T$, $f$ factors through
a projective object $P$ in $\mathcal{A}$. Therefore the morphism
$(0, f)$ factors through $i_1(P)$,  which is projective in ${\rm
Mon}(\mathcal{A})$; see Lemma \ref{lem:Frobeniusexact}. Hence
$(0,f)=0$ in
 the stable monomorphism category $\underline{\rm
Mon}(\mathcal{A})$.
\vskip 3pt

To check (T2) for $T'$, note that each object $\alpha$ fits into a
conflation $$i_2(s(\alpha))\longrightarrow \alpha \longrightarrow i_1({\rm
Cok}\; \alpha)$$
 and thus into a triangle
 $$i_2(s(\alpha))\longrightarrow \alpha
\longrightarrow i_1({\rm Cok}\; \alpha) \longrightarrow
i_2(s(\alpha))[1].$$
Here as in section 2, $s(\alpha)$ denotes the
source of $\alpha$. Hence the smallest triangulated subcategory of
$\underline{\rm Mon}(\mathcal{A})$ containing
$i_1(\underline{\mathcal{A}})$ and $i_2(\underline{\mathcal{A}})$ is
$\underline{\rm Mon}(\mathcal{A})$ itself. Now applying the
condition (T2) of $T$, we infer that (T2) holds for $T'$.

\vskip 3pt

 Finally to
see the condition (T3) for $T'$, it is direct to check that ${\rm
End}_{\underline{\rm Mon}(\mathcal{A})}(T')\simeq T_2({\rm
End}_{\underline{\mathcal{A}}}(T))$. Note that the algebra ${\rm
End}_{\underline{\mathcal{A}}}(T)$ has finite global dimension. Then
by \cite[Chapter III, Proposition 2.6]{ARS} we infer that ${\rm
End}_{\underline{\rm Mon}(\mathcal{A})}(T')$ has finite global
dimension.
\end{proof}

We will give an application of Theorem \ref{thm:first}.
Let $A=\oplus_{n\geq 0}A_n$ be a positively graded artin algebra.
Denote by $c$ the maximal integer such that $A_c\neq 0$. Consider the following upper triangular
matrix algebra $$\mathrm{b}(A)=\begin{pmatrix} A_0 & A_1 & \cdots & A_{c-2} & A_{c-1}\\
                     0     &  A_0  & \cdots   & A_{c-3} & A_{c-2}\\
                       \vdots &  \vdots &   \ddots & \vdots &  \vdots \\
                       0 & 0  & \cdots  & A_0 & A_1\\
                         0 & 0 & \cdots  & 0 & A_0 \end{pmatrix}. $$
Here the multiplication of $\mathrm{b}(A)$ is induced from the one
of $A$. This algebra is called the \emph{Beilinson algebra} of $A$ in \cite{C2}.

\vskip 5pt

Denote by ${\rm mod}^\mathbb{Z}\; A$ the category of finitely
generated $\mathbb{Z}$-graded $A$-modules with homomorphisms
preserving degrees. We say that $A$ is \emph{graded self-injective}
provided that ${\rm mod}^\mathbb{Z}\; A$ is a Frobenius category. In
fact, this is equivalent to that as a ungraded algebra $A$ is
self-injective (\cite{GG,C2}). In this case, we denote by
$\underline{\rm mod}^\mathbb{Z}\; A$  the stable category of ${\rm
mod}^\mathbb{Z}\; A$ modulo projective modules; it is a triangulated
category.

\vskip 5pt

 We say that  a graded algebra $A$ is  \emph{right
well-graded}, provided that $A_c$, as a right $A_0$-module, is
sincere in the sense of \cite[p.317]{ARS}. In fact, for a graded
self-injective algebra $A$, it is right well-graded if and only if
it is left well-graded; see \cite[Lemma 2.2]{C2}. In this case we
will simply say that the graded algebra $A$ is \emph{well-graded}.

\begin{cor}\label{cor:first}
Let $A=\oplus_{n\geq 0} A_n$ be a positively graded self-injective artin algebra which is well-graded. Suppose that
$A_0$ has finite global dimension. Then there is a triangle equivalence
$$\underline{\rm Mon}({\rm mod}^\mathbb{Z} \; A)\; \simeq \; \mathbf{D}^b({\rm mod}\; T_2(\mathrm{b}(A))).$$.
\end{cor}

\begin{proof}
Note that by \cite[Chapter III, Proposition 2.6]{ARS}, the Beilinson
algebra $\mathrm{b}(A)$ and then $T_2(\mathrm{b}(A))$ has finite
global dimension. By \cite[Corollary 1.2]{C2} there is a triangle
equivalence
 $\underline{{\rm mod}}^\mathbb{Z}\; A \simeq \mathbf{D}^b({\rm mod}\; \mathrm{b}(A))$.
 In particular, there is a tilting object $T$ in $\underline{{\rm mod}}^\mathbb{Z}\; A$ with
endomorphism algebra $\mathrm{b}(A)$. We apply  Theorem
\ref{thm:first} to get a tilting object $T'$ in $\underline{\rm
Mon}({\rm mod}^\mathbb{Z}\; A)$  whose endomorphism algebra is
isomorphic to $T_2(\mathrm{b}(A))$. Note that the stable
monomorphism category $\underline{\rm Mon}({\rm mod}^\mathbb{Z} \;
A)$ is idempotent-split; in fact, it is even a Krull-Schmidt
category. Then the result follows immediately from Lemma
\ref{lem:Keller}.
\end{proof}

\section{Stable Monomorphism Category as  Singularity Category}

In this section we will relate the stable monomorphism category of
the (graded) module category of a (graded) self-injective algebra to
the (graded) singularity category of the associated (graded)
$2\times 2$ upper triangular matrix algebra.

\vskip 5pt

Let $A$ be an artin algebra. Recall that the bounded homotopy
category $\mathbf{K}^b({\rm proj}\; A)$ of projective modules is
viewed naturally as a triangulated subcategory of $\mathbf{D}^b({\rm
mod}\; A)$.  Following \cite{O1,O2}, we call the Verdier quotient
triangulated category
$$\mathbf{D}_{\rm
sg}(A)=\mathbf{D}^b({\rm mod}\; A)/\mathbf{K}^b({\rm proj}\; A)$$
the \emph{singularity category} of $A$; compare \cite{Buc} and
\cite{H2}.

\vskip 5pt

 Our second main observation is as follows. Recall that for an artin algebra $A$, $T_2(A)$
 is the $2\times 2$ upper triangular matrix algebra of $A$.

\begin{thm}\label{thm:second}
Let $A$ be a self-injective algebra.  Then we have a
triangle equivalence
$$\underline{\rm Mon}({\rm mod}\; A)\; \simeq \; \mathbf{D}_{\rm sg}(T_2(A)).$$
\end{thm}

Before giving the proof, we recall several notions. Let $A$ be an
artin algebra. Following \cite[p.400]{AM} an acyclic complex
$P^\bullet$ of projective $A$-modules is called \emph{totally
acyclic} if the Hom complex ${\rm Hom}_A(P^\bullet, A)$ is still
acyclic (also see \cite[section 7]{Kr}). An $A$-module $M$ is said
to be \emph{Gorenstein-projective} if there is a totally acyclic
complex $P^\bullet$ such that its zeroth cocycle $Z^0(P^\bullet)$ is
isomorphic to $M$ (\cite[Chapter 10]{EJ}). Note that a module $M$ is
Gorenstein-projective if and only if ${\rm Ext}^i_A(M, A)=0$, ${\rm
Ext}_{A^{\rm op}}^i({\rm Hom}_A(M, A), A)=0$ for $i\geq 1$ and the
natural map $M \rightarrow {\rm Hom}_{A^{\rm op}}({\rm Hom}_A(M, A),
A)$ is an isomorphism (compare \cite[Definition (1.1.2)]{Chr}).
Denote by $\mbox{Gproj\; }A$ the full subcategory of $\mbox{mod}\;
A$ consisting of Gorenstein-projective modules. Note that projective
modules are Gorenstein-projective and thus $\mbox{proj}\; A
\subseteq \mbox{Gproj}\; A$. Moreover, by \cite[Proposition 5.1]{AR}
the subcategory ${\rm Gproj}\;A$ is closed under extensions and
taking direct summands (also see \cite{EJ}); and then it is direct
to see that ${\rm Gproj}\; A$ is a Frobenius exact category such
that its projective objects are equal to projective $A$-modules
(\cite[Proposition 3.1(1)]{C3}). Denote by $\underline{\rm Gproj}\;
A$ its stable category modulo projective $A$-modules; it is a
triangulated category.

\vskip 5pt

Recall that  an artin algebra $A$ is said to be \emph{Gorenstein} if
the regular modules $_AA$ and $A_A$  have finite injective
dimensions (\cite{H2}). In this case the two dimensions are equal
and the common value is denoted by ${\rm G.dim}\; A$. We say that
the Gorenstein algebra $A$ is \emph{$1$-Gorenstein} provided that
${\rm G.dim}\; A\leq 1$.

\vskip 5pt

For an artin algebra $A$, denote by ${\rm sub}\;A$ the full
subcategory of ${\rm mod}\;A$ consisting of submodules of projective
modules; these modules are called \emph{torsionless modules}.

\vskip 5pt

 We note the following well-known result.

\begin{lem}\label{lem:wellknown}
Let $A$ be a $1$-Gorenstein algebra. Then we have ${\rm Gproj}\;
A={\rm sub}\;A$.
\end{lem}

\begin{proof}
The inclusion ${\rm Gproj}\; A\subseteq {\rm sub}\;A$ is easy. On
the other hand, assume that $M$ is  a torsionless module. Consider a
short exact sequence $0\rightarrow M\rightarrow P\rightarrow
M'\rightarrow 0$ with $P$ projective.  Since the regular module
$A_A$ has injective dimension $1$, using dimension shift, we infer
that ${\rm Ext}^i(M, A)=0$ for $i\geq 1$. Then by \cite[Corollary
11.5.3]{EJ} (see also \cite[Lemma 3.7]{C3} and \cite[Proposition
7.13]{Kr}), $M$ is Gorenstein-projective.
\end{proof}

The next observation is essentially due to Li and Zhang
(\cite[Theorem 1.1]{LZ}; also see \cite[Proposition 3.6]{IKM}). Recall that for an artin algebra $A$, a
(right) module over $T_2(A)$ is identified with a morphism of
(right) $A$-modules; in fact, this yields an equivalence ${\rm
mod}\; T_2(A)\simeq {\rm Mor}({\rm mod}\; A)$ of categories; see
\cite[Chapter III, Proposition 2.2]{ARS}.

\begin{lem}\label{lem:LiZhang}
Let $A$ be a self-injective algebra. Then we have an equivalence  of
exact categories
$${\rm Mon}({\rm mod}\; A) \; \simeq  \; {\rm Gproj}\; T_2(A).$$
\end{lem}

\begin{proof}
 We claim  first that under the equivalence ${\rm
mod}\; T_2(A)\simeq {\rm Mor}({\rm mod}\; A)$, ${\rm sub}\; T_2(A)$
corresponds to ${\rm Mon}({\rm mod}\; A)$. Note that the regular
module $T_2(A)_{T_2(A)}$ corresponds to the monomorphism $\binom{0}
{{\rm Id}_A}\colon A\rightarrow A\oplus A$. From this one infers
that torsionless $T_2(A)$-modules correspond to monomorphisms in
${\rm mod}\; A$. On the other hand, the third paragraph of the proof
of Lemma \ref{lem:Frobeniusexact} already shows that for a
monomorphism $\alpha$, there is a short exact sequence $0\rightarrow
\alpha \rightarrow \binom{0}{{\rm Id}_P} \rightarrow
\alpha'\rightarrow 0$ in ${\rm Mor}({\rm mod}\; A)$ such that $P$ is
a projective $A$-module. Note that $\binom{0}{{\rm Id}_P}$
corresponds to a projective $T_2(A)$-module. Therefore the
monomorphism $\alpha$ corresponds to a torsionless $T_2(A)$-module.
This completes the proof of the claim.

\vskip 3pt

 Note that by \cite[Remark 3.5]{C1} (also see \cite{FGR,
H2}) the algebra $T_2(A)$ is $1$-Gorenstein and then we can apply
Lemma \ref{lem:wellknown}. By the claim above, we obtain an
equivalence ${\rm Mon}({\rm mod}\; A) \simeq
 {\rm Gproj}\; T_2(A)$ of categories. It is direct to check that this equivalence preserves the exact
 structures, that is, the equivalence and its inverse preserve
 short exact sequences in ${\rm Mon}({\rm mod}\; A) $ and ${\rm Gproj}\; T_2(A)$.
\end{proof}

We will recall the last ingredient in our proof. Let $A$ be an artin
algebra. Consider the following composite of natural functors
$$F_A\colon \mbox{Gproj}\; A \hookrightarrow \mbox{mod}\; A\hookrightarrow \mathbf{D}^b(\mbox{mod}\; A)
\longrightarrow \mathbf{D}_{\rm sg}(A)$$ where from the left side,
the first functor is the inclusion, the second identifies modules
with stalk complexes concentrated in degree zero (\cite[p.40,
Proposition 4.3]{Har}) and the last is the quotient functor. Note
that the additive functor $F_A$ vanishes on projective modules and
then induces uniquely an additive  functor
$\mbox{\underline{Gproj}}\; A \rightarrow \mathbf{D}_{\rm sg}(A)$,
which is still denoted by $F_A$. The following important result is
due to Buchweitz (\cite[Theorem 4.4.1]{Buc}) and independently to
Happel (\cite[Theorem 4.6]{H2}); also see \cite[Proposition 3.5 and
Theorem 3.8]{C3}.

\begin{lem} \label{lem:BuchweitzHappel} {\rm (Buchweitz-Happel)} Let $A$
be an artin algebra. Then the functor $F_A\colon
\mbox{\underline{\rm Gproj}}\; A\rightarrow \mathbf{D}_{\rm sg}(A)$
is a fully faithful triangle functor. Moreover, if $A$ is
Gorenstein, then the functor $F_A$ is dense and thus a triangle
equivalence.
\end{lem}

\vskip 5pt

\noindent {\bf Proof of Theorem \ref{thm:second}.}\quad By Lemma
\ref{lem:LiZhang} there is an equivalence ${\rm Mon}({\rm mod}\; A)
\simeq {\rm Gproj}\; T_2(A)$ of Frobenius exact categories.  Hence
we have an induced equivalence of triangulated categories
$$\underline{\rm
Mon}({\rm mod}\; A)\; \simeq \; \underline{\mbox{{\rm Gproj}}}\;T_2(A). $$
Note that $T_2(A)$ is Gorenstein by \cite[Theorem 3.3]{C1} (also see \cite{FGR, H2}).
Thus the result follows
directly from Lemma \ref{lem:BuchweitzHappel}. \hfill $\square$

\vskip 7pt

We will also need a graded version of Theorem \ref{thm:second}. Let
$A=\oplus_{n\geq 0}A_n$ be a positively graded artin algebra. Denote by
 ${\rm proj}^\mathbb{Z}\; A$ the full
subcategory of ${\rm mod}^\mathbb{Z}\; A$ consisting of projective objects.
Following \cite{O2} one has the \emph{graded
singularity category} of $A$ defined by
$$\mathbf{D}_{\rm
sg}^\mathbb{Z}(A) =\mathbf{D}^b({\rm mod}^\mathbb{Z}\;
A)/\mathbf{K}^b({\rm proj}^\mathbb{Z}\; A).$$

\vskip 5pt

 For a graded  module $M=\oplus_{i\in \mathbb{Z}} M_i$ and
an integer $d\in \mathbb{Z}$, its \emph{shifted module} $M(d)$ has
the same module structure as $M$ while it is graded such that
$M(d)_i=M_{d+i}$ for all $i\in \mathbb{Z}$. This defines
automorphisms $(d)\colon {\rm mod}^\mathbb{Z}\; A \rightarrow {\rm
mod}^\mathbb{Z}\; A$, which are called \emph{degree-shift functors}.
For graded modules $M, N$, we write ${\rm HOM}_A(M, N)=\oplus_{i\in
\mathbb{Z}}{\rm Hom}_{{\rm mod}^\mathbb{Z}\; A}(M, N(i))$ and set
${\rm EXT}_A^n(-, -)$ to be the $n$-th right derived functors
(\cite{NV} and \cite{GG}).

\vskip 5pt

 An acyclic complex $P^\bullet$ in ${\rm
proj}^\mathbb{Z}\; A$ is \emph{totally acyclic} if the Hom complex
${\rm HOM}_A(P^\bullet, A)$ in ${\rm proj}^\mathbb{Z}\; A^{\rm op}$
is still acyclic. A graded $A$-module is called  \emph{graded
Gorenstein-projective} provided that it is the zeroth cocycle of a
totally acyclic complex. Thus we have a full subcategory ${\rm
Gproj}^\mathbb{Z}\; A$ of ${\rm mod}^\mathbb{Z}\; A$ consisting of
graded Gorenstein-projective modules and evidently ${\rm
proj}^\mathbb{Z}\; A\subseteq {\rm Gproj}^\mathbb{Z}\; A$. As in the
ungraded case, the category ${\rm Gproj}^\mathbb{Z}\; A$ is a
Frobenius exact category with its projective objects
 equal to graded projective $A$-modules.

\vskip 5pt

 Recall that a graded artin algebra $A$ is said to be \emph{graded Gorenstein} if
the graded regular modules $_AA$ and $A_A$ have finite injective
dimensions in ${\rm mod}^\mathbb{Z}\; A$ and ${\rm mod}^\mathbb{Z}\;
A^{\rm op}$, respectively. In this case the two dimensions are the
same, which will be denoted by ${\rm G.dim}^\mathbb{Z}\; A$.

\vskip 5pt

 We note the following fact, which guarantees in principle that most results in  Gorenstein
 homological algebra hold true in the graded situation.

\begin{lem}
Let $A$ be a positively graded artin algebra, and let $M$ be a
graded $A$-module. Then we have
\begin{enumerate}
\item[(1)] the module $M$ is graded Gorenstein-projective if and
only if it is Gorenstein-projective as a ungraded module;
\item[(2)] the algebra $A$ is graded Gorenstein if and  only if it
is Gorenstein as a ungraded algebra; in this case, we have ${\rm
G.dim}^\mathbb{Z}\; A={\rm G.dim}\; A$.
\end{enumerate}
\end{lem}

\begin{proof}
For (1), it suffices to note that a graded module $M$ is graded
Gorenstein-projective if and only if  ${\rm EXT}^i_A(M, A)=0$, ${\rm
EXT}_{A^{\rm op}}^i({\rm HOM}_A(M, A), A)=0$ for $i\geq 1$ and the
natural map $M \rightarrow {\rm HOM}_{A^{\rm op}}({\rm HOM}_A(M, A),
A)$ is an isomorphism of graded modules; moreover, for graded
modules $M$ and $N$ we have for each $i$ a natural identification
${\rm EXT}^i_A(M, N)={\rm Ext}^i_A(M, N)$ (\cite[Corollary
2.4.7]{NV}). For (2), note that a graded module $M$ has finite
injective dimension in ${\rm mod}^\mathbb{Z}\; A$ if and only if it
has finite injective dimension as a ungraded module; moreover, the
two dimensions are the same (\cite[Theorem 2.8.7]{NV}).
\end{proof}

One can show the  graded analogues of Lemmas \ref{lem:wellknown},
\ref{lem:LiZhang} and \ref{lem:BuchweitzHappel}. Using these, we have  the
following graded analogue of Theorem \ref{thm:second}.

\begin{prop}\label{prop:first} Let $A=\oplus_{n\geq 0} A_n$ be a positively graded self-injective artin
 algebra. Denote by $T_2(A)$ the $2\times 2$ upper triangular matrix algebra of $A$ which is
graded such that $T_2(A)_n=T_2(A_n)$ for $n\geq 0$. Then we have a triangle
equivalence
$$\underline{\rm Mon}({\rm mod}^\mathbb{Z}\; A)\; \simeq \;\mathbf{D}^\mathbb{Z}_{\rm sg}(T_2(A)).$$
\end{prop}

\section{The Stable Category of Ringel-Schmidmerier}

In the last section, we will give the promised two characterizations to  the stable
category of Ringel-Schmidmeier (\cite{RS3}).
\vskip 5pt

Let $k$ be a field and let $p\geq 2$ be an integer. In particular, the fixed
commutative artinian ring $R$ now is the field $k$. Consider the truncated polynomial algebra
 $A=k[t]/(t^p)$ with $t$ an indeterminant; it is positively graded such that ${\rm deg}\; t=1$. Note that
 $A$ is graded self-injective and moreover it is well-graded.
 In particular, the category ${\rm mod}^\mathbb{Z}\; A$ of finitely generated
 graded $A$-modules is Frobenius. Following \cite[subsection 0.4]{RS3}, we denote by
 $\mathcal{S}(\widetilde{p})$ the category of pairs $(V, U)$, where $V$
is a graded module over $A$ and $U\subseteq V$ is a graded
submodule, and the morphisms in this category are given by morphisms
in the graded module category which  respect the inclusion. There is
a natural identification $\mathcal{S}(\widetilde{p})={\rm Mon}({\rm
mod}^\mathbb{Z}\; A)$ and then by Lemma \ref{lem:Frobeniusexact} it
is a Frobenius exact category.  Hence its stable category
$\underline{\mathcal{S}}(\widetilde{p})$ modulo projective objects
is triangulated. It will be called the \emph{stable category of
Ringel-Schmidmeier}.

\vskip 5pt

We note that the Beilinson algebra $\mathrm{b}(A)$ (in the sense of \cite{C2})
of  the graded algebra $A$ is isomorphic to the path algebra $k\mathbb{A}_{p-1}$ of the linear quiver $\mathbb{A}_{p-1}$ with $p-1$
vertices and linear orientation (compare \cite[Example 2.9]{O2}).
 Then the $2\times 2$ upper triangular matrix algebra
$\Lambda=T_2(\mathrm{b}(A))$ is given by the following quiver with
$2p-2$ vertices subject to the commutativity  relation
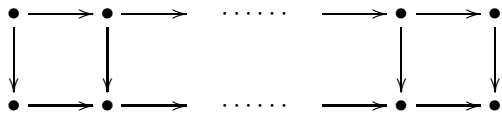
\[\xymatrix{
\bullet \ar[r] \ar[d] & \bullet \ar[d]\ar[r] & \quad \cdots \cdots
\quad \ar[r] &  \bullet  \ar[r] \ar[d]&
\bullet \ar[d] \\
\bullet \ar[r] & \bullet \ar[r] & \quad \cdots \cdots \quad \ar[r] &
\bullet \ar[r] & \bullet}
\]
Let us mention that similar diagrams appeared  in \cite{Lad}.

\vskip 5pt

Recall that the $2\times 2$ upper triangular matrix algebra $T_2(A)$
is graded such that $T_2(A)_n=T_2(A_n)$ for $n\geq 0$. Then it is
isomorphic, as a graded algebra, to  $T_2(k)[t]/(t^p)$, while the
latter is graded such that ${\rm deg} \; T_2(k)=0$ and ${\rm deg}\;
t=1$.

\vskip 5pt

Combining Corollary \ref{cor:first} and Proposition \ref{prop:first}
together, we obtain two characterizations of the stable category of
Ringel-Schmidmeier.

\begin{thm}\label{thm:last}
Use the notation above. Then there are triangle equivalences
$$\mathbf{D}^b({\rm mod}\; \Lambda)\;  \simeq \; \underline{\mathcal{S}}(\widetilde{p})
 \; \simeq \; \mathbf{D}^\mathbb{Z}_{\rm sg}(T_2(k)[t]/(t^p)).$$
\end{thm}

\vskip 5pt

\begin{rem}
Let us remark that taking into account of the results obtained in \cite{KMLP} and \cite[Corollary 1.2]{Lad},
one may find a close relation between some results in \cite{KML} and Theorem \ref{thm:last}.
\hfill $\square$
\end{rem}

\vskip 10pt

\noindent {\bf Acknowledgements} \quad The author would like to
thank Dirk Kussin very much for a private communication. In fact,
this work is inspired by a talk given by Dirk Kussin in Paderborn.
The author is indebted to Henning Krause for  illuminating
discussions and to Bernhard Keller and Zhi-Wei Li for helpful
comments.

\bibliography{}

\begin{thebibliography}{999}

\bibitem{Ar} {\sc D.M. Arnold,} Abelian Groups and Representations of Finite Partially Ordered
Sets, Springer GMS Books in Mathematics, 2000.

\bibitem{AR} {\sc M. Auslander and I. Reiten,} {\em Applications of contravariantly finite
subcategories}, Adv. Math. {\bf 86} (1)(1991), 111--152.


\bibitem{ARS}{\sc M. Auslander, I. Reiten and S.O. Smal{\o},}
Representation Theory of Artin Algebras, Cambridge Univ. Press,
1995.

\bibitem{AM} {\sc L.L. Avramov and A. Martsinkovsky,} {\em Absolute, relative, and Tate cohomology of
modules of finite Gorenstein dimension,} Proc. London Math. Soc. (3)
{\bf 85} (2002), no.2, 393--440


\bibitem{B} {\sc G. Birkhoff}, {\em Subgroups of abelian groups}, Proc. London Math.
Soc. II, Ser. {\bf 38} (1934), 385--401.

\bibitem{BK} {\sc A.I. Bondal and M.M. Kapranov,} {\em Enhanced triangulated categories}, Math.
USSR Sbornik {\bf 70} (1990), no.1, 93--107.


\bibitem{Buc} {\sc  R.O. Buchweitz},  Maximal Cohen-Macaulay Modules and
Tate Cohomology over Gorenstein Rings, Unpublished manuscript, 155pp
(1987). Available at https://tspace.library.utoronto.ca/handle/1807/16682.


\bibitem{C1} {\sc X.W. Chen}, {\em Singularity categories, Schur functors and triangular matrix
rings}, Algebra and Representation Theory {\bf 12} (2009), 181--191.



\bibitem{C2} {\sc X.W. Chen}, {\em Graded self-injective algebras ``are" trivial
extensions,} J. Algebra {\bf 322} (2009), 2601--2606.



\bibitem{C3} {\sc X.W. Chen}, {\em Relative singularity categories and Gorenstein-projective modules},
Math. Nach., to appear, arXiv: 0709.1762.


\bibitem{Chr} {\sc L.W. Christensen}, Gorenstein Dimension, Lecture Notes in Math. {\bf 1747}, Springer-Verlag,
Berlin Heidelberg, 2000.


\bibitem{GG} {\sc R. Gordon and E.L. Green,}
{\em Graded Artin algebras,} J. Algebra {\bf 76}(1) (1982),
111--137.


\bibitem{EJ} {\sc E.E. Enochs and O.M.G. Jenda}, Relative Homological Algebra, de
Gruyter Expositions in Math. {\bf 30}, Walter de Gruyter, Berlin New
York, 2000.


\bibitem{FGR} {\sc R.M. Fossum, P. Griffith  and I. Reiten}, Trivial Extensions of Abelian Categories, Lecture
Notes in Math. {\bf 456}, Springer-Verlag, Berlin Heidelberg New York, 1975.



\bibitem{H1}{\sc D. Happel,} Triangulated Categories in the
Representation Theory of Finite Dimensional Algebras, London
Mathematical Society Lecture Note Ser. {\bf 119}, Cambridge Univ.
Press, 1988.

\bibitem{H2} {\sc D. Happel},  {\em On Gorenstein algebras,}
Progress in Math., vol. {\bf 95}, 389--404, Birkh$\ddot{a}$user
Verlag, Basel, 1991.



\bibitem{Har} {\sc R. Hartshorne}, Residue and Duality,
Lecture Notes in Math. {\bf 20}, Springer-Verlag, 1966.

\bibitem{IKM} {\sc O. Iyama, K. Kato and J.I. Miyachi,} {\em
Recollement of homotopy categories and Cohen-Macaylay modules},
arXiv:0911.0172v1.


\bibitem{Ke3} {\sc B. Keller},  {\em Chain complexes and stable
categories}, Manuscripta Math. {\bf 67} (1990), 379--417.


\bibitem{Ke3.5} {\sc B. Keller}, {\em Derived categories and their uses},
in: Handbook of Algebra {\bf 1}, 671--701, North-Holland, Amsterdam,
1996.



\bibitem{Ke4} {\sc  B. Keller}, {\em Derived categories and
tilting}, in: Handbook of Tilting Theory, 49--104,
London Math. Soc. Lecture Note Ser. {\bf 332},
Cambridge Univ. Press, Cambridge, 2007.



\bibitem{Kr}  {\sc H. Krause}, {\em The stable derived category of a noetherian scheme},
Compositio Math. {\bf 141} (2005), 1128--1162.



\bibitem{KML} {\sc D. Kussin, H. Lenzing and H. Meltzer}, {\em Categories of vector bundles and invariant subspaces of nilpotent operators}, preprint.


\bibitem{KMLP} {\sc D. Kussin, H. Lenzing, H. Meltzer and  J.A. de la Pe\~na}, {\em Stable categories of
vector bundles}, preprint.



\bibitem{Lad} {\sc S. Ladkani,} {\em On derived equivalences of lines, rectangles and triangles,} 
arXiv:0911.5137.



\bibitem{Len} {\sc H. Lenzing}, {\em Hereditary categories}, in: Handbook of Tilting Theory, 105--146,
London Math. Soc. Lecture Note Ser. {\bf 332},
Cambridge Univ. Press, Cambridge, 2007.


\bibitem{LZ} {\sc Z.W. Li and P. Zhang}, {\em A construction of Gorenstein-projective modules,}
submitted for publication, 2009.




\bibitem{NV}{\sc C. Nastasescu and F. Van Oystaeyen,} Methods of Graded
Rings, Lecture Note in Math. {\bf 1836}, Springer, 2004.



\bibitem{O1} {\sc D. Orlov}, {\em Triangulated categories of
singularities and D-branes in Landau-Ginzburg models}, Proc. Steklov
Inst. Math. {\bf 246} (3) (2004), 227--248.



\bibitem{O2}{\sc D. Orlov}, {\em Derived categories of coherent
sheaves and triangulated categories of singularities,}
math.AG/0503632v2.


\bibitem{RS1} {\sc C.M. Ringel and M. Schmidmeier}, {\em Submodule categories of wild representation type},
Jour. Pure and Applied Algebra {\bf 205} (2) (2006), 412--422.


\bibitem{RS2} {\sc C.M. Ringel and M. Schmidmeier}, {\em The Auslander-Reiten translation in submodule categories},
Trans. Amer. Math. Soc. {\bf 360} (2008), 691--716.

\bibitem{RS3} {\sc C.M. Ringel and M. Schmidmeier}, {\em Invariant subspaces of nilpotent operators I},
Journal Reine Angew. Math. Band  {\bf 614} (2008), 1--52.

\bibitem{Sim} {\sc D. Simson,} {\em Chain categories of modules and subprojective representations of posets
over uniserial algebras}, Rocky Mountain J. Math. {\bf 32} (2002), 1627--1650.


\end{thebibliography}

\vskip 10pt

 {\footnotesize \noindent Xiao-Wu Chen, Department of
Mathematics, University of Science and Technology of
China, Hefei 230026, P. R. China \\
Homepage: http://mail.ustc.edu.cn/$^\sim$xwchen \\
\emph{Current
address}: Institut fuer Mathematik, Universitaet Paderborn, 33095,
Paderborn, Germany}

\end{document}